\documentclass[12pt]{amsart}
\usepackage[dvips]{epsfig}
\usepackage{graphics}
\usepackage{latexsym}
\usepackage{verbatim}
\usepackage{amsmath}
\usepackage{amsthm}
\usepackage{amssymb}
\newtheorem{theorem}{Theorem}[section]
\newtheorem{lemma}[theorem]{Lemma}

\newtheorem{proposition}[theorem]{Proposition}

\newcommand{\ens}[1]{\mathbb{#1}}

\newcommand{\R}{\mathbb{R}}

\def\cal{\mathcal}

\def\derpar#1#2{\frac{\partial#1}{\partial#2}}
\def\var{\varepsilon}
\def\signcm{\bigskip\bigskip\hspace{80mm}
\vbox{{\sc C. Mouhot\par\vspace{3mm}
CEREMADE, Univ. Paris IX Dauphine \par
Place du Mar\'echal de Lattres de Tassigny \par
75775 Paris cedex 16 \par
FRANCE \par\vspace{3mm}
e-mail:} cmouhot@ceremade.dauphine.fr }}
\def\signrs{\bigskip\bigskip\hspace{80mm}
\vbox{{\sc R. M. Strain\par\vspace{3mm}
Harvard University\par
Department of Mathematics\par
One Oxford Street\par
Cambridge, MA 02138\par
USA\par\vspace{3mm}
e-mail:} strain@math.harvard.edu}}

\begin{document}

\title[Coercivity estimates for non-cutoff Boltzmann operators]
{Spectral gap and coercivity estimates for linearized Boltzmann collision operators without angular cutoff}

\author{Cl\'ement Mouhot \& Robert M. Strain}
\date{\today}

\hyphenation{bounda-ry rea-so-na-ble be-ha-vior pro-per-ties
cha-rac-te-ris-tic}

\begin{abstract}
In this paper we prove new constructive coercivity estimates for the Boltzmann  
collision operator without cutoff, that is for {\em long-range} interactions. 
In particular we give a generalized sufficient condition for the existence of a spectral 
gap which involves both the growth behavior of the collision kernel at large relative velocities and 
its singular behavior at grazing and frontal collisions. It provides in particular existence of a spectral gap 
and estimates on it for interactions deriving from the hard potentials $\phi(r)=r^{-(s-1)}$, $s \ge 5$ 
or the so-called moderately soft potentials $\phi(r)=r^{-(s-1)}$, $3< s < 5$, (without angular cutoff). 
In particular this paper recovers (by constructive means), improves and extends previous 
results of Pao~\cite{Pao74}. We also obtain constructive coercivity estimates for 
the Landau collision operator for the optimal coercivity norm pointed out in~\cite{Guo:Land} 
and we formulate a conjecture about a unified necessary and sufficient condition 
for the existence of a spectral gap for Boltzmann and Landau linearized collision operators. 
\end{abstract}

\maketitle



\tableofcontents

\section{Introduction}
\setcounter{equation}{0}

This paper deals with the question of spectral gap and coercivity estimates for the 
Boltzmann and Landau integro-differential collision operators. This is motivated by 
the question of obtaining explicit constant in the recent new energy methods 
in~\cite{Guo:VPB,Guo:Land,Guo:VMB,Guo:BEsoft,Guo:whole,StGu:rel,StGu:almostexp,MoNe,StGu:expsoft,Strain}. 
The starting points are the recent new constructive tools of works~\cite{BaMo,Mcoerc}.  
As we shall explain this work recovers, improves, makes explicit and clarifies the works 
of Pao~\cite{Pao74} three decades ago. 
Before entering into the details, let us introduce the mathematical objects. 

\subsection{The Boltzmann equation}

The Boltzmann equation (Cf. \cite{Ce88,CIP}) 
describes the behavior of a dilute gas when the only 
interactions taken into account are binary collisions. It reads 
in some space domain $\Omega \subset \R^N$:
 \begin{equation}\label{e1}
 \derpar{f}{t} + v \cdot \nabla_x f = Q_{\cal{B}} (f,f), \qquad  x \in \Omega, \quad v \in \R^N, \quad t \geq 0,
 \end{equation}
where $N \ge 2$ is the dimension. 
In equation~\eqref{e1}, $Q$ is the quadratic
{\it Boltzmann collision operator}, defined by
$$
Q_{\mathcal{B}}(f,f)= \int_{\R^N \times \ens{S}^{N-1}} \big[ f(v') \, f(v'_*) \,-f(v) \, f(v_*) \,\big] \, 
B(|v-v_*|, \cos \theta) \, dv_* \, d\sigma
$$
in the so-called ``$\sigma$-representation" (see~\cite[Chapter~1, Section~4.6]{Vi:hb}. 
In this representation the parametrization of the collision is 
 \[ v' = \frac{v+v_*}2 + \frac{|v-v_*|}2 \sigma, \quad v'_* = \frac{v+v_*}2 - \frac{|v-v_*|}2 \sigma, 
     \qquad \sigma \in \ens{S}^{N-1} \]
and the {\it deviation angle} is defined by $\cos \theta = (v'_*-v') \cdot (v_*-v) /|v_*-v|^2$. 
Up to a jacobian factor $2^{N-2} \, \sin^{N-2} (\theta/2)$ (see again~\cite[Chapter~1, Section~4.6]{Vi:hb}),
one can also define the alternative so-called ``$\omega$-representation"
 \[ Q_{\mathcal{B}}(f,f)= \int_{\R^N \times \ens{S}^{N-1}} \big[ f(v') \, f(v'_*) \,-f(v) \, f(v_*) \,\big] 
       \, 2^{N-2} \, \sin^{N-2} (\theta/2) \, B(|v-v_*|, \cos \theta) \, dv_* \, d\omega \]
with the formula
 \[ v' = v+ ((v_*-v) \cdot \omega) \, \omega, \quad v'_* = v_* - ((v_*-v)\cdot \omega) \, \omega, 
    \qquad \omega \in \ens{S}^{N-1}. \]
Remark that this operator is local in $x$, and therefore all its functional study in the sequel shall be 
done with no space variable $x$.   
\smallskip

Boltzmann's collision operator has the fundamental properties of
conserving mass, momentum and energy
 \begin{equation*}
 \int_{\R^N}Q_{\cal{B}}(f,f) \, \phi(v)\,dv = 0, \qquad
 \phi(v)=1,v,|v|^2
 \end{equation*}
and satisfies well-known Boltzmann's $H$ theorem
 \begin{equation*}
 - \frac{d}{dt} \int_{\R^N} f \log f \, dv = - \int_{\R^N} Q(f,f)\log(f) \, dv \geq 0.
 \end{equation*}
The functional $- \int f \log f$ is the {\em entropy} of the
solution. Boltzmann's $H$ theorem implies that at some point 
$x \in \Omega$, any equilibrium
distribution function, {\it i.e.}, any function which is a maximum of the
entropy, has the form of a locally Maxwellian distribution
 \begin{equation}
\label{maxw}
 M(\rho,u,T)(v)=\frac{\rho}{(2\pi T)^{N/2}}
 \exp \left( - \frac{\vert u - v \vert^2} {2T} \right), 
 \end{equation}
where $\rho,\,u,\,T$ are the {\em density}, {\em mean velocity}
and {\em temperature} of the gas at the point $x$, defined by
 \begin{equation}
\label{field}
 \rho = \int_{\R^N}f(v) \, dv, \quad u =
 \frac{1}{\rho}\int_{\R^N} v \, f(v) \, dv, \quad T = {1\over{N \rho}}
 \int_{\R^N}\vert u - v \vert^2 \, f(v) \, dv.
 \end{equation}
For further details on the physical background and derivation of
the Boltzmann equation we refer to~\cite{Ce88,CIP} and~\cite{Vi:hb}.
\smallskip

We consider collision kernels of the form 
\begin{equation}\label{tenscond}
B(|v-v_*|, \cos \theta) = |v-v_*|^\gamma \, b(\cos \theta), \quad \gamma \in (-N,+\infty)
\end{equation}
with 
\begin{equation}\label{singcond}
b(\cos \theta) \sim_{\theta \sim 0}  b^*(\theta) \, \left(\sin \theta/2 \right)^{-(N-1)-\alpha}, 
   \quad \alpha \in [0,2), 
\end{equation}
where $b^*(\theta)$ is non-negative, bounded and non-zero near $\theta \sim 0$. 
When $\alpha \ge 0$ the angular singularity is not integrable, the operator 
is said to be {\em non-cutoff}. 

An important remark on the quadratic collision operator is: by using the 
change of variable $\sigma \to - \sigma$ (which changes $\theta$ into 
$\pi - \theta$) 
one can replace $b$ by 
$$
\tilde b(\cos \theta) = \frac12 \, {\bf 1}_{\theta \in [0,\pi/2]} \, 
\left[ b(\cos \theta) + b (\cos ( \pi - \theta)) \right]
$$ 
(where ${\bf 1}_E$ denotes the usual characteristic function of a set $E$. 
Therefore we shall consider without restriction in the sequel that $b$ 
satisfies the singularity condition~\eqref{singcond} above, 
and is $0$ on $[\pi/2,\pi]$. 
\smallskip

For particles interacting according to some spherical repulsive potential   
$$
\phi(r)=r^{-(s-1)} , \quad s \in [2,+\infty),
$$
the collision kernel is not explicit but it can be shown that for the dimension $N=3$,  
$B$ satisfies~\eqref{tenscond}  with $\gamma = (s-5)/(s-1)$ 
and~\eqref{singcond}  with  $\alpha = 2/(s-1)$ (see for instance~\cite{Ce88,CIP,ADVW}). 
Therefore as a convention one shall speak in~\eqref{tenscond} of {\em hard potentials}   
when $s \ge 5$, {\em maxwellian potential}  when $s=5$, 
{\em soft potentials}   when $2<s <5$. Moreover among soft potentials we shall denote by  
{\em moderately soft potentials}  the case when $3 \le s <0$. 

Let us mention also for the sake of completeness that in the case of contact interactions 
(the so-called {\it hard spheres} model), the collision kernel 
is locally integrable and explicit: it takes the form (in dimension $N=3$) $B(|v-v_*|, \cos \theta) = |v-v_*|$ (up to 
a normalization constant). 

\subsection{Linearized Boltzmann collision operator}

We denote by 
$$
\mu =\mu(v) :=(2\pi)^{-N/2}e^{-|v|^2/2}
$$ 
the normalized unique equilibrium 
with mass $1$, momentum $0$ and temperature $1$. 
We consider fluctuations around this equilibrium of the form 
$$
f = \mu + \mu^{1/2} \, g
$$
which results the following linearized collision operator (note the sign convention): 
 \[ L_{\cal{B}}(g) = - \mu^{-1/2} \, \left[ Q_{\cal{B}}(\mu,\mu^{1/2} g) + Q_{\cal{B}}(\mu^{1/2} g, \mu) \right]. \]
 
For the sake of simplification we shall always consider in the sequel the linearized collision 
operator around this normalized equilibrium. This is no restriction: a detailed 
discussion of the dependence of the spectral gap and coercivity estimates on this operator 
in terms of the mass, mean velocity and temperature can be found in~\cite{Mcoerc}. 

It is well-known (see~\cite{Ce88} for instance) that $L_{\cal{B}}$ (acting in the velocity space) 
is an unbounded symmetric operator on $L^2$, such that its Dirichlet form satisfies 
$$
D_{\cal{B}} (g) := \langle L_{\cal{B}}g, g\rangle \ge 0,
$$ 
and that $D_{\cal{B}} (g) = 0$ if and only if $g={\bf P}g$ where
$$
{\bf P}g=\left( a+b\cdot v+c|v|^2\right) \mu^{1/2}
$$
(with $a,c\in\mathbb{R}$ and $b\in \mathbb{R}^N$) is the $L^2$ orthogonal 
projection onto the space of the so-called ``collisional invariants"
$$
\mbox{Span} \big\{ \mu^{1/2}, v_1 \mu^{1/2}, \dots, v_N \mu^{1/2}, |v|^2 \mu^{1/2} \big\}.
$$

\subsection{The Landau equation} 

The Landau equation was written by Landau in 1936 (\cite{Land36}). It is similar to the 
Boltzmann equation~\eqref{e1} but with a different collision operator $Q_{\cal{L}}$. 
Indeed in the case of long-distance interactions, collisions occur 
mostly for very small $\theta$. When all collisions become 
concentrated on $\theta =0$, one obtains by the so-called 
{\em grazing collision limit} asymptotic 
(see for instance \cite{ArBu:91,DeLu:92,Desv:asBE:92,Vill:nocu:98,AlVi:04} for a detailed discussion) 
the {\em Landau collision operator} 
\begin{equation}
Q_{\mathcal{L}}(f,f) =\nabla \cdot \left(\int_{{\mathbb R}^N}{\bf a}(v-v_*)[f_* (\nabla f)-f(\nabla f)_*]dv_*\right),  \label{lan} 
\end{equation}
where $\nabla=(\partial_{v_1},\ldots,\partial_{v_N})$, $f=f(v)$, $f_*=f(v_*)$ etc. The non-negative symmetric $N\times N$ matrix ${\bf a}= {\bf a}(z)$ is 
\begin{equation}
a_{ij}(z)=\left\{ \delta _{ij}-\frac{z_iz_j}{|z|^2}\right\} |z|^{\gamma
+2}.  \label{phi}
\end{equation}

This operator is used for instance in plasma physics in the case of a Coulomb 
potential where $\Phi(|z|) = |z|^{-3}$ in dimension $N=3$ 
(for more details see~\cite[Chapter~1, Section~1.7]{Vi:hb} and the references therein). 
Indeed let us mention that for Coulomb interactions the Boltzmann collision operator 
does not make sense anymore (see~\cite[Annex~I, Appendix]{Vi:habil}). 
By analogy with the Boltzmann's case, and even if the physical meaning 
of such collision kernel is not clear for the Landau collision operator, 
one shall speak of {\em hard potentials}  when $\gamma > 0$, 
{\em maxwellian potentials}  when $\gamma =0$, and {\em soft potentials} 
when $\gamma < 0$. Moreover among soft potentials we shall denote by 
{\em moderately soft potentials}  the case when $-2 \le \gamma <0$. 

\subsection{Linearized Landau collision operator}
Consider again fluctuation around the equilibrium of the form
$$
f=\mu + \mu^{1/2} \, g.
$$
Then the linearized Landau collision operator is defined by  
$$
L_{\mathcal{L}}g = - \mu^{-1/2}  \left[ Q_{\cal{L}}(\mu,\mu^{1/2} g) + Q_{\cal{L}}(\mu^{1/2} g, \mu)\right].
$$

It was proved in \cite{Guo:Land,Hinton,DeLe97} that that $L_{\cal{L}}$ (acting on the velocity space) 
is an unbounded symmetric operator on $L^2$, such that its Dirichlet form satisfies 
$$
D_{\cal{L}} (g) := \langle L_{\cal{L}}g, g\rangle \ge 0
$$ 
and $D_{\cal{L}} (g) = 0$ if and only if $g={\bf P}g$.

\subsection{Previous results and motivations}
The study of spectral gap estimates for the linearized 
Boltzmann and Landau collision operators has a long 
history, see for instance~\cite{Hilb:EB:12,WCUh52,Carl:57,Grad58,Grad63,WCUh70,Cafl80,Boby88} 
and we refer to~\cite{BaMo,Mcoerc} for a more detailed discussion of it. 
Let us just emphasize some references directly related to this paper. For linearized 
Boltzmann collision operators with locally integrable collision kernel (which is  satisfied 
for instance under the so-called {\em Grad's angular cutoff}), the existence of a 
spectral gap is equivalent to $\gamma \ge 0$ (see~\cite{Grad63,Cafl80,CIP} for 
non-constructive proofs). The question of obtaining polynomial or exponential rate 
of relaxation of the form $O(e^{-\lambda t^\beta})$ for $0<\beta<1$ and for soft potentials with cutoff was studied in~\cite{Cafl80,UkAs,StGu:almostexp,StGu:expsoft}, and 
explicit spectral gap and generalized coercivity estimates were obtained in~\cite{BaMo,Mcoerc}. 

For the sake of completeness let us also mention that Cercignani introduced 
an alternative ``potential cutoff" to Grad's angular cutoff in~\cite{Ce67} for which results 
on the spectrum of the linearized collision operator were obtained in~\cite{CzPa}. 

While the cutoff theory of the linearized collision operator (without space variable) is quite developed, 
for non locally integrable collision kernels, there are quite few works. First for the Boltzmann collision 
operator at the linearized level there
are essentially the two papers by Pao~\cite{Pao74}: he proved that the resolvent is compact for 
inverse power-laws interaction potentials with $s > 3$ (in our notation). His proof was using tools 
from  pseudo-differential operators theory: the idea was to compute the symbol of the linearized  
operator $L_{\cal{B}}$, and then to search for adequate condition for $L_{\cal{B}} + C$ (for some 
constant $C >0$) to be invertible with compact inverse, by reducing to the Maxwell case and 
decomposing along spherical harmonics. This fundamental work was critically reviewed 
in particular for the use of techniques from pseudo-differential operator theory in~\cite{Klau77}.

The critics of~\cite{Klau77} on~\cite{Pao74}, together with the fact that the proof in~\cite{Pao74} 
was highly technical, were probably the reasons for which the results in~\cite{Pao74} were not 
considered as completely reliable, and this paper was somehow forgotten in the following decades. 

One of our goal in this work is to clarify the discussion about the validity of the results~\cite{Pao74}, 
by recovering and improving strongly these results. Moreover we replace the question treated in~\cite{Pao74} 
into the unified framework of quantitative coercivity estimates for the linearized collision operators, 
which is related to the works~\cite{BaMo,Mcoerc} (for the Landau collision operator 
at the linearized level, let us mention the key works~\cite{DeLe97,Guo:Land} and the approach 
developed in~\cite{Mcoerc} which shall be used here). A general motivation for this framework is to 
provide explicit rate of convergence to equilibrium in the energy methods 
which have emerged recently in the collisional kinetic theory: see~\cite{Guo:VPB,Guo:Land,Guo:VMB,Guo:BEsoft,Guo:whole,StGu:rel,StGu:almostexp,MoNe,StGu:expsoft,Strain}. 

At the non-linear level, let us mention some breakthroughs related to non-cutoff interactions: 
in the spatially homogeneous Cauchy theory:~\cite{ArkInf1,ArkInf2,ElmInf,Vill:nocu:98,GoNewClass}, 
in the grazing collision limit 
from Boltzmann to Landau equation:~\cite{ArBu:91,Desv:asBE:92,DeLu:92}, in the parabolic-like 
regularizing property in the velocity variable:
~\cite{Desv:reg:95,Desv:reg:96,Desv:reg:97,Lions98,Vi:reg:98,ADVW}, in the construction of renormalized 
solutions (in the spirit of~\cite{DiLiBE}):~\cite{AlViCPAM,AlVi:04}.   

In the spirit of the method of Pao, a systematic approach by Fourier transform and pseudo-differential 
operators for the study of linear and non-linear Boltzmann non-cutoff collision operators has been 
developed by Alexandre~\cite{AlexTTSP,AlexM2AN,AlexM3AS}. The applications of this technical 
tools are not clear at now, but they seem to have contributed to the understanding of the Boltzmann 
equation in two ways: in the development of a renormalization process adapted to the non-cutoff case
(\cite{AlexCRASrenorm,AlViCPAM}), and by providing some sharp estimates from above 
on the collision operator. 
 
\subsection{Main theorems} 
First we state our coercivity results for the linearized Boltzmann collision operator: 
\begin{theorem}\label{theo:mainB}
Let $B$ be a collision kernel satisfying~(\ref{tenscond},\ref{singcond}). Then 
\begin{itemize}
\item For any $\var > 0$ there is a constant $C_{B,\var}$, constructive from our proof and 
depending on $B$ and $\var$, such that the Dirichlet form $D_{\cal{B}}$ of the linearized 
Boltzmann collision operator associated to $B$ satisfies 
$$
D_{\cal{B}}(g) \ge C_{B,\var} \, \big\| [ g - {\bf P} g ] \, (1+|v|^2)^{(\gamma+\alpha-\var)/4} \big\|^2 _{L^2(\R^N)}.
$$
\item There is a constant $C_{B,0}$ (obtained by non-constructive means in our proof) such that 
$$
D_{\cal{B}}(g) \ge C_{B,0} \, \big\| [ g - {\bf P} g ] \, (1+|v|^2)^{(\gamma+\alpha)/4} \big\|^2 _{L^2(\R^N)}.
$$
\end{itemize}
\end{theorem}

Second we state the constructive version of the coercivity result in~\cite{Guo:Land} for the 
linearized Landau collision operator: 
\begin{theorem}\label{theo:mainL}
Let $\gamma \in [-N,+\infty)$. Then there is some constant $C_\gamma$, constructive from our proof and depending on $\gamma$, such that the  Dirichlet form $D_{\cal{L}}$ of the associated linearized 
Landau collision operator satisfies: 
$$
D_{\cal{L}}(g) \ge C_\gamma \, \big\| [ g - {\bf P}g ] \big\|_\sigma ^2.  
$$
where $\|\cdot \|_\sigma$ is the following anisotropic norm:  
\begin{multline*}
\|g\|_{\sigma}^2 := \Big\|  (1+|v|^2)^{\gamma/4} \, \Pi_v\nabla_v g \Big\|^2 _{L^2(\R^N)} \\
+ \Big\| (1+|v|^2)^{(\gamma +2)/4} \, [I-\Pi_v] \nabla_v g \Big\|^2 _{L^2(\R^N)} 
+\Big\| (1+|v|^2)^{(\gamma+2)/4} \, g \Big\|^2 _{L^2(\R^N)} 
\end{multline*}
with
$$
\Pi_v \nabla_v g= \left(\frac{v}{|v|}\cdot \nabla_v g\right) \frac{v}{|v|}.
$$
\end{theorem}

\subsection{Consequences on the spectrum}
In a previous paper the first author has proved: 
\begin{theorem}[\cite{Mcoerc}] \label{theo:reg}
Let $B$ be a collision kernel satisfying~(\ref{tenscond},\ref{singcond}) with $\alpha >0$. Then 
there is a constant $C_{B}$ (constructive from the proof and depending on $B$) 
such that the Dirichlet form $D_{\cal{B}}$ of the linearized Boltzmann collision operator associated to $B$ 
satisfies 
$$
D_{\cal{B}}(g) \ge C_{B} \, \big\| [ g - {\bf P}g ] \big\|^2 _{H^{\alpha/2} _{\mbox{\scriptsize {\em loc}}} (\R^N)}.
$$
\end{theorem}

Therefore as soon as $\gamma + \alpha > 0$ and $\alpha >0$ it is straightforward by gathering 
Theorems~\ref{theo:mainB} and~\ref{theo:reg} that the resolvent of $L_{\cal{B}}$ is compact 
and so that the spectrum is purely discrete and the eigenvectors basis is complete in $L^2$. 
In dimension $N=3$ for interactions deriving from inverse-power laws potentials $\phi(r)=r^{-(s-1)}$, 
the corresponding condition on $s$ is: $(s-5)/(s-1) + 2/(s-1) >0$, that is $s > 3$. 
Hence we can recover completely 
the results of Pao~\cite{Pao74}, by constructive means and without tools of pseudo-differential operators. 

Moreover Theorem~\ref{theo:mainB}  yields the first estimates on this spectral gap in the case 
$\gamma + \alpha >0$ and answers to the question of the existence of a spectral gap 
in the limit case $\gamma + \alpha =0$. 
\smallskip

For the linearized Landau collision operator, it had been already shown in~\cite{Guo:Land} that 
$$
C_1 \, \big\| [ g - {\bf P}g ] \big\|_\sigma ^2 \le D_{\cal{L}}(g) \le C_2 \, \big\| [ g - {\bf P}g ]\big\|_\sigma ^2 
$$
by non-constructive means (in fact the constant $C_2$ could be made explicit from the proof 
whereas the constant $C_1$ was obtained by a compactness argument). 
This implies straightforwardly that the resolvent is compact as 
soon as $\gamma >-2$ (which implies in this case that the spectrum is purely discrete and 
the eigenvectors basis is complete in $L^2$). Moreover spectral gap exists if and only 
if $\gamma \ge -2$. 

In this context, Theorem~\ref{theo:mainL} provides the first estimates on it for  
this whole region (as well as an explicit version of the coercivity estimate). 

\subsection{Conjecture and perspectives} 
For the linearized Boltzmann collision operator with angular cutoff, the existence of a spectral gap 
is equivalent to $\gamma \ge 0$ (see~\cite{Grad63,Cafl80,CIP}). 
This situation can be loosely thought of as part of the limit case ``$\alpha =0$". For the linearized 
Landau collision operator, the existence of a spectral gap is equivalent to $\gamma \ge -2$ 
as discussed above. This situation can be thought of as the limit case ``$\alpha =2$". 

From the necessary and sufficient condition in these two limit cases, and the sufficient condition 
of Theorem~\ref{theo:mainB} in the intermediate cases, it is natural to state the following conjecture: 
\smallskip

\noindent 
{\bf Conjecture.} {\em Let $B$ be a collision kernel satisfying~(\ref{tenscond},\ref{singcond}) 
with $\gamma \in (-N, +\infty)$ and $\alpha \in [0,2)$. 
Then the linearized Boltzmann collision operator associated to $B$ admits 
a spectral gap if and only if $\gamma + \alpha \ge 0$. Moreover this statement is still valid if one 
includes formally the case of angular cutoff in ``$\alpha =0$", and add the linearized Landau collision 
operator as the limit case ``$\alpha =2$".} 
\smallskip

Let us remark that in the limit cases of the linearized Boltzmann collision operator with angular cutoff and the 
linearized Landau collision operator the conjecture is proved 
as discussed before. And Theorem~\ref{theo:mainB}  proves that the condition $\gamma + \alpha \ge 0$ 
is sufficient for the existence of a spectral gap.  

Here are now some open questions linked with this conjecture: 
\begin{enumerate}
\item In order to show that the condition is necessary, is it possible to find some particular sequence 
of functions $g_n$ contradicting the spectral gap estimate when $\gamma + \alpha <0$? 
\item A related and more difficult question is to identify and understand the coercivity norm for 
the non-cutoff linearized 
Boltzmann collision operator. This norm is likely to be intricate and anisotropic for the weight on the 
fractional derivative part (as for the linearized Landau collision operator). This point 
is related to the open question of constructing smooth perturbative solutions near equilibrium for the 
Boltzmann equation without angular cutoff (see the discussion in~\cite{MoNe}). 
\item The spectrum is purely discrete in the case $\gamma + \alpha >0$, $\alpha >0$ (and we conjecture that this 
is also true on the non-cutoff borderline case $\gamma + \alpha >0$, $\alpha =0$). In the angular cutoff 
case the geometry of the spectrum can be obtained by perturbation arguments (Cf.~\cite{Grad63,Cafl80,CIP}). 
Hence there remains the region $\gamma + \alpha \le 0$ for which there is no information on the 
geometry of the spectrum at now (including the linearized Landau collision operator when $\alpha=2$ 
and $\gamma < -2$).  
\end{enumerate}

Finally we shall put in perspective this conjecture with the non-linear case. A non-linear 
analogous of a spectral gap is provided by Cercignani's conjecture: 
$$
\cal{D}_{\cal{B}} (f) \ge C_f \, H(f | M_f) 
$$
where 
$$
H(f|M_f) = H(f) - H(M_f) := \int_{\R^N} f \, \log \frac{f}{M_f} \, dv
$$ 
is the relative entropy between $f$ and its associated Maxwellian equilibrium distribution, 
and 
$$
\cal{D}_{\cal{B}} (f) := - \int_{\R^N} Q_{\cal{B}}(f,f) \, \log f \, dv \quad ( \ge 0) 
$$ 
denotes the entropy production functional (for more details we refer to~\cite{Vi03}). 
It was shown in~\cite{DeVi:L2:00} that an equivalent version of this conjecture is true for the 
Landau collision operator when $\gamma \ge 0$, and it was shown in~\cite{Vi03}  that 
this conjecture for the Boltzmann collision operator for 
$B(v-v_*,\cos \theta) \ge K_B \, (1+ |v-v_*|)^\gamma$ with $\gamma \ge 2$. From this two limit 
cases Villani ``interpolated" the following conjecture: 
\smallskip

{\bf Conjecture (Villani~\cite{Vi03}).} {\it Cercignani's conjecture is satisfied if and only if 
$\gamma + \alpha \ge 2$. }
\smallskip

(The other part of the conjecture on the question of the existence of a spectral
gap for the linearized collision operator was not correct --as shown by Theorem~\ref{theo:mainB}--,
which is likely to be explained by the confusion in the field about Pao's results). 
Obviously this conjecture leaves more room than the one we made on the spectral gap, since 
it is not stated in which functional space the distribution $f$ lives,  
but from the results in ~\cite{Vi03} and~\cite{BobyCerc:conjCer:99}, 
it is reasonable to try to prove that Villani's conjecture holds in  
spaces of nonnegative functions in $L^1(1+|v|^q) \cap H^k$ for any $q,k \ge 0$. 

Still an interesting 
question is to know whether Cercignani's conjecture could still be true for $\gamma + \alpha <2$ 
in some functional space with ``stretched exponential decay" (see the discussions 
in~\cite{BobyCerc:conjCer:99,Mrate}). 
This is important since the condition $\gamma + \alpha \ge 2$ rules out all physical interactions. 
 
Another important question, if Villani's and our conjectures hold, would be to understand the 
reasons for this ``gap" between the non-linear and linearized behavior of the entropy production 
for the Boltzmann equation. 

\subsection{Methods of the proof and plan of the paper} 
Section~\ref{sec:tech} is devoted to a technical estimate of decay on the kernel 
of the non-local part of the linearized Boltzmann collision operator, 
under the angular cutoff assumption. We show that its mixing effects 
can be quantified into a gain of polynomial weights (this phenomenon 
was noticed in~\cite{Grad63} and is here fully developed). 
This estimate has been isolated from the rest of the proof since it is the key step, 
and also since it can be of independent interest for researchers in the field. 

Section~\ref{sec:Bconst} is devoted to the (constructive) proof of the coercivity estimate for the linearized 
Boltzmann collision operator when $\var >0$ (first point in Theorem~\ref{theo:mainB}). 
The idea is to estimate from below the Dirichlet form by truncating the angular part 
$b$ of the collision kernel on the angles $\theta \in [0,\theta_0(|v-v_*|)]$ with 
$\theta_0(|v-v_*|) \sim |v-v_*|^{-k}$ for some suitable $k>0$, and then 
balance the lower bound on the local part (in which the polynomial weight 
results from a competing effect between the fact the lower bound on $b$ is 
big in this region, and the fact that the size of this angular region is small), and the upper 
bound on the non-local part, for which the previous technical estimate plays a crucial 
role. The constructive coercivity estimates from~\cite{Mcoerc} are also used. 

Section~\ref{sec:Bnc} is devoted to the (non-constructive) proof of the coercivity 
estimate for the linearized Boltzmann collision operator when $\var =0$ (second point in 
Theorem~\ref{theo:mainB}). The idea is to reduce to a cutoff-like linearized Boltzmann collision 
operator with a different collision kernel, and then apply a strategy based on  
Weyl's Theorem about compact perturbation of essential spectrum, in the spirit of~\cite{Grad63,CIP}.

Finally Section~\ref{sec:L} is devoted to the (constructive) proof of the coercivity estimate 
for the linearized Landau collision operator (Theorem~\ref{theo:mainL}). The idea is to combine 
the estimates in~\cite{Guo:Land} on the different parts of the decomposition of the operator 
as ``diffusive part + bounded part" (which involves the optimal coercivity norms), 
with the constructive coercivity estimates of~\cite{Mcoerc} obtained on the global operator, 
but with non-optimal coercivity norms. 

\section{A technical estimate on $L_{\cal{B}}$}\label{sec:tech}
\setcounter{equation}{0}

We assume in this section that the collision kernel $B$ takes the particular form 
(sometimes named ``variable hard spheres collision kernels") 
$$
B_q (|v-v_*|, \cos \theta) = |v-v_*|^q, \quad q \in (-N,+\infty). 
$$
These collision kernels are non-physical except in the case $q =1$ (hard spheres) 
but they shall play an important role in intermediate steps of our proof in the next sections. They 
satisfy the angular cutoff assumption, that is $B_q$ is locally integrable 
in terms of $v,v_*,\sigma$. It is well-known that under this assumption the collision 
operator can be split into ``gain" and ``loss" parts. 

Therefore one can decompose the linearized collision operator 
$L_{\cal{B}}$ corresponding to $B_q$ as follows 
$$
L_{\cal{B}}g =  \nu_{\cal{B}} \, g - K_{\cal{B}}g
$$
with the non-local part 
\begin{multline*}
K_{\cal{B}}g := \int_{\R^N \times \ens{S}^{N-1}} \Big[ g(v') \, \mu^{1/2} (v'_*) \, 
+ g(v'_*) \, \mu^{1/2} (v') - \\ g(v_*) \, \mu^{1/2}(v) \Big] \, \mu^{1/2}(v_*) \,  
B_q(|v-v_*|) \, dv_* \, d\sigma
\end{multline*}
and the following multiplicative function for the local part 
$$
\nu_{\cal{B}}(v) := \int_{\R^N \times \ens{S}^{N-1}} \mu(v) \, B_q(|v-v_*|) \, dv_* \, d\sigma 
= \big|\ens{S}^{N-1}\big| \, \big( |\cdot|^q * \mu \big) (v).
$$ 
The non-local part $K_{\cal{B}}$ itself splits into a pure convolution part 
$$
K_{\cal{B}} ^c g := \big|\ens{S}^{N-1}\big| \, \big[ |\cdot|^q * (\mu^{1/2} \, g) \big] (v) \, \mu^{1/2}(v)
$$
and a ``gain" part (denoted so since it flows from the gain part of the original non-linear operator) 
$$
K_{\cal{B}} ^+ g := \int_{\R^N \times \ens{S}^{N-1}} \Big[ g(v') \, \mu^{1/2} (v'_*) \, 
+ g(v'_*) \, \mu^{1/2} (v') \Big] \, \mu^{1/2}(v_*) \,  B_q(|v-v_*|) \, dv_* \, d\sigma.
$$
Note that using the change of variable $\sigma \to -\sigma$, the latter writes also 
$$
K_{\cal{B}} ^+ g := 2 \, \int_{\R^N \times \ens{S}^{N-1}} g(v') \, \mu^{1/2} (v'_*) \, 
\mu^{1/2}(v_*) \,  B_q(|v-v_*|) \, dv_* \, d\sigma.
$$
\smallskip

We shall focus on the kernel of the $K_{\cal{B}} ^+$ part. 
For any locally integrable collision kernel $B$ we can define 
$k_B=k_B(v,v')$ such that 
$$ 
\forall \, v \in \R^N, \quad K_{\cal{B}} ^+g(v) = \int_{\R^N} g(v') \, k_B (v,v') \, dv'
$$
and we denote by $k_q:=k_{B_q}$ the kernel obtained for the particular collision kernel 
$B_q$ above. Note that this definition is independent of the representation 
we use for the collision operator, either using the ``$\sigma$" or ``$\omega$" 
parametrization of the unit vector on the sphere. 
\medskip

Then we have the
 \begin{proposition}\label{prop:Grad}
 The kernel $k_q$ is symmetric:
 $$
 \forall \, v,v' \in \R^N, \quad k_q(v,v') = k_q(v',v),
 $$
 and for $q>-1$ and for any $s \in \R$ it satisfies the control 
 $$
  \forall \, v \in \R^N, \quad \int_{\R^N} k_q(v,v') \, (1+|v'|)^s \, dv' \le
         C_{q,s} \, (1+|v|)^{q+s-(N-1)}
 $$
 for some explicit constant $C_{q,s}$ depending only on $q,s$.
 \end{proposition}

Let us first compute the kernel $k_q$ using well-known changes of variable 
(see for instance~\cite{Grad63}).

 \begin{lemma}\label{lem:kernel}
 We have, for $q>-1$, the explicit formula (let us recall that $\omega=(v'-v)/|v'-v|$):
   \begin{multline*}
   k_q(v,v') = \frac{2^N}{|v'-v| \, (2 \pi)^{N/2}} \,
             \exp\left\{ -\frac{|v'-v|^2}8 -\frac{|v'-v + 2 (v \cdot \omega) \omega|^2}8 \right\} \\ \times
                 \left( \int_{\omega^\bot} |v'-v + z|^{q-(N-2)}
                 \exp\left\{ - \frac{|z + \left(v - (v \cdot \omega) \omega \right)|^2}2 \right\} \, dz \right).
   \end{multline*}
 \end{lemma}

\begin{proof}[Proof of Lemma~\ref{lem:kernel}]
We start from the $\omega$-representation of $K_{\cal{B}} ^+$: 
 \begin{multline*}
 I:= 2 \, \int_{\R^N \times \ens{S}^{N-1}}
     g(v') \, (\mu_* \mu_* ')^{1/2} \, B_q \, dv_* \, d\sigma \\
   = 2^{N-1} \, \int_{\R^N \times \ens{S}^{N-1}}
     g(v') \, (\mu_* \mu_* ')^{1/2} \, |v-v_*|^q \, \sin^{N-2} (\theta/2) \, dv_* \, d\omega \\
   = 2^{N-1} \, \int_{\R^N \times \ens{S}^{N-1}}
     g(v') \, (\mu_* \mu_* ')^{1/2} \, |v-v_*|^{q-(N-2)} \, |v'-v|^{N-2} \, dv_* \, d\omega.
 \end{multline*}

Then we perform the change variable $v_* \to V=v_*-v$ (with $\omega$ fixed), and
the change of variable $V = r \, \omega + z$ with $z \in \omega^\bot$ (with $\omega$ fixed).
These two changes of variables have jacobian equal to $1$. We obtain
 \[ I =  2^{N-1} \, \int_{\R \times \ens{S}^{N-1}}
     g(v+r \,\omega) \, |r|^{N-2} \,
     \left( \int_{\omega^\bot} (\mu_* \mu_* ')^{1/2} \, 
     |r \omega + z|^{q-(N-2)} \, dz \right) \, dr \, d\omega. \]
Finally we make the spherical change of variable $(r,\omega) \in \R \times \ens{S}^{N-1}
\to W=r\omega \in \R^N$ with jacobian $2 \, r^{-(N-1)}$ (a factor $2$ comes from the fact
that $r \in \R$), which yields
 \[ I =  2^{N} \, \int_{\R^N} g(v+W) \, |W|^{-1} \,
     \left( \int_{W^\bot} (\mu_* \mu_* ')^{1/2} \, |W+z|^{q-(N-2)} \, dz \right) \, dW. \]

Let us rewrite the argument of the exponential term:
 \begin{multline*}
 |v_*|^2 + |v' _*|^2 = |v+W+z|^2 + |v+z|^2 = \frac12 |W + 2(v+z)|^2 + \frac12 |W|^2 \\
     = \frac12 |W + 2 (v \cdot \omega) \omega|^2
    + 2 |z + \left(v - (v \cdot \omega) \omega \right)|^2 + \frac12 |W|^2
 \end{multline*}
where $\omega = W/|W|$ (in the last equality we have used that $z \, \bot \, W$).
We deduce that
 \[ (\mu_* \mu_* ')^{1/2} = \frac{1}{(2 \pi)^{N/2}} \, 
       \exp\left\{ -\frac{|W|^2}8 - \frac{|z + \left(v - (v \cdot \omega) \omega \right)|^2}2
                            -\frac{|W + 2 (v \cdot \omega) \omega|^2}8 \right\}. \]
We deduce the formula. One checks that the integral on $z$ is well-defined
for $q>-1$ since in this case $q-(N-2) >-(N-1)$.
\end{proof}

\begin{proof}[Proof of Proposition~\ref{prop:Grad}]
The fact that $k_q$ is symmetric is easily checked from the formula. 
Let us turn to the bound from above. 
\smallskip

First let us assume that $q-(N-2) \ge 0$. 
We estimate the integral over $z \in \omega^\bot$ 
in terms of $v$:
 \begin{multline*}
 \int_{\omega^\bot} |v'-v + z|^{q-(N-2)}
                 \exp\left\{ -\frac{|z + \left(v - (v \cdot \omega) \omega \right)|^2}2 \right\} \, dz \\
 = \int_{\omega^\bot} |v'-v + \bar z - \left(v - (v \cdot \omega) \omega \right)|^{q-(N-2)}
                 \exp\left\{ -\frac{|\bar z|^2}2 \right\} \, d\bar z \\
 \le C \, (1+|v - (v \cdot \omega) \omega |)^{q-(N-2)} \, (1+|v'-v|)^{q-(N-2)} \\
 \le C \, (1+|v|)^{q-(N-2)} \, (1+|v'-v|)^{q-(N-2)}.
 \end{multline*}

We deduce that
 \begin{multline*}
 \int_{\R^N} k_q(v,v') \, (1+|v'|)^s \, dv' \\
 \le C \, (1+|v|)^{q-(N-2)} \,
     \int_{\R^N} |v'-v|^{-1} \, (1+|v'|)^s \,  (1+|v'-v|)^{q-(N-2)} \\ 
             \exp\left\{ -\frac{|v'-v|^2}8 -\frac{|v'-v + 2 (v \cdot \omega) \omega|^2}8 \right\} \, dv'
 \end{multline*}
and we use the elementary inequality
 \[ (1+|v'|)^s \le C \, (1+|v|)^s \, (1+|v'-v|)^{|s|} \]
to get
 \[  \int_{\R^N} k_q(v,v') \, (1+|v'|)^s \, dv' \le C \, (1+|v|)^{q+s-(N-2)} \, J \]
with
 \begin{multline*}
 J := \int_{\R^N} |v'-v|^{-1} \, (1+|v'-v|)^{|s| + q-(N-2)} \\
             \exp\left\{ -\frac{|v'-v|^2}8 -\frac{|v'-v + 2 (v \cdot \omega) \omega|^2}8 \right\} \, dv'.
 \end{multline*}

Now in the term $J$ we perform the change of variable $v' \to u=v'-v$ and then 
the spherical change of variables $u = r \, \omega$, $r \in \R_+$, $\omega \in \ens{S}^{N-1}$, 
choosing $v$ as the north pole vector in the angle parametrization (note that this is not 
the same parametrization as the one used to define the deviation angle). It yields
 \[ J = |\ens{S}^{N-2}| \, \int_0 ^{+\infty} r^{N-2} \, (1+r)^{|s|+q-(N-2)} \, e^{-\frac{r^2}8}
        \int_0 ^\pi e^{-\frac{(r+2|v| \cos \varphi)^2}8} \, \sin^{N-2} \varphi \, d\varphi \, dr. \]

Then for technical reasons we treat separately the case $N=2$ and $N \ge 3$.
\smallskip

First for $N \ge 3$, we have $\sin^{N-2} \varphi \le \sin \varphi$, and thus
 \[ J \le C \, \int_0 ^{+\infty} r^{N-2} \, (1+r)^{|s|+q-(N-2)} \, e^{-\frac{r^2}8}
        \int_0 ^\pi e^{-\frac{(r+2|v| \cos \varphi)^2}8} \, \sin \varphi \, d\varphi \, dr. \]
Then we do the change of variable $y=(r+2|v|\cos \varphi)$ in the $\varphi$ integral, to
get
 \[ J \le C \, |v|^{-1} \, \int_0 ^{+\infty} r^{N-2} \, (1+r)^{|s|+q-(N-2)} \, e^{-\frac{r^2}8} \, dr
        \int_{-\infty} ^{+\infty} e^{-\frac{y^2}8} \, dy \le C \, |v|^{-1} \]
which yields the conclusion for large $v$ (the estimate for small $v$ is
immediate).
\smallskip

Second for $N=2$,  we perform the same changes of variable which yields
 \begin{multline*}
 J \le C \, |v|^{-1} \, \int_0 ^{+\infty} \, (1+r)^{|s|+q-(N-2)} \, e^{-\frac{r^2}8} \\
        \int_{(r-2|v|)} ^{(r+2|v|)} e^{-\frac{y^2}8} \,
        \left( 1 - \left(\frac{y -r}{2 \, |v|}\right)^2 \right)^{-1/2} \, dy \, dr \\
 \le C \, \int_0 ^{+\infty} \, (1+r)^{|s|+q-(N-2)} \, e^{-\frac{r^2}8}
        \int_{(r-2|v|)} ^{(r+2|v|)} e^{-\frac{y^2}8} \, \left( 4 \, |v|^2 - (y -r)^2 \right)^{-1/2} \, dy \, dr.
 \end{multline*}
Finally we split into two parts: $|y -r| \le |v|$ and $|y -r| \ge |v|$.
On the first part we have $\left( 4 \, |v|^2 - (y -r)^2 \right)^{-1/2} \le |v|^{-1}$
which yields the result. On the second part we have either $r \ge |v|/2$
or $|y| \ge |v|/2$ which gives an exponential decay in $v$ thanks to the terms 
$e^{-\frac{r^2}8}$ and $e^{-\frac{y^2}8}$ in the integrand. This concludes the proof
in dimension $N=2$.
\smallskip

Now let us come back to the case $q-(N-2) <0$. Then we have
 \begin{multline*}
 \int_{\omega^\bot} |v'-v + z|^{q-(N-2)}
                 \exp\left\{ -\frac{|z + \left(v - (v \cdot \omega) \omega \right)|^2}2 \right\} \, dz \\
 = \int_{\omega^\bot} |v'-v + \bar z - \left(v - (v \cdot \omega) \omega \right)|^{q-(N-2)}
                 \exp\left\{ -\frac{|\bar z|^2}2 \right\} \, d\bar z \\
 \le C \, (1+|v- (v \cdot \omega) \omega|)^{q-(N-2)} \, (1+|v'-v|)^{|q-(N-2)|}.
 \end{multline*}
The additional difficulty will be therefore to obtain decay in $v$ since the weight in this formula 
only involves the projection of $v$ on $\omega^\bot$.

Let us follow the computations as before: we have
 \begin{multline*}
 \int_{\R^N} k_q(v,v') \, (1+|v'|)^s \, dv' \\ \le C \,
     \int_{\R^N} |v'-v|^{-1} \, (1+|v'|)^s \, 
      (1+|v- (v \cdot \omega) \omega|)^{q-(N-2)} \, (1+|v'-v|)^{|q-(N-2)|} \\
            \times  \exp\left\{ -\frac{|v'-v|^2}8 -\frac{|v'-v + 2 (v \cdot \omega) \omega|^2}8 \right\} \, dv'.
 \end{multline*}
We use again the inequality
 \[ (1+|v'|)^s \le C \, (1+|v|)^s \, (1+|v'-v|)^{|s|} \]
to get
 \[  \int_{\R^N} k_q(v,v') \, (1+|v'|)^s \, dv' \le C \, (1+|v|)^s \, J' \]
with
 \begin{multline*}
 J' := \int_{\R^N} |v'-v|^{-1} \, (1+|v'-v|)^{|s|+|q-(N-2)|} \, 
         (1+|v- (v \cdot \omega) \omega|)^{q-(N-2)} \\
             \exp\left\{ -\frac{|v'-v|^2}8 -\frac{|v'-v + 2 (v \cdot \omega) \omega|^2}8 \right\} \, dv'.
 \end{multline*}

Again we perform the change of variable $v' \to u=v'-v$ and then the spherical 
change of variables $u = r \, \omega$, $r \in \R_+$, $\omega \in \ens{S}^{N-1}$, choosing $v$ as the
north pole vector in the angle parametrization. It yields
 \begin{multline*} 
 J' = |\ens{S}^{N-2}| \, \int_0 ^{+\infty} r^{N-2} \, (1+r)^{|s|+|q-(N-2)|} \, e^{-\frac{r^2}8} \\
   \int_0 ^\pi (1+|v| \, \sin \varphi)^{q-(N-2)} \,
      e^{-\frac{(r+2|v| \cos \varphi)^2}8} \, \sin^{N-2} \theta \, d\varphi \, dr. 
 \end{multline*}
Then we split between $|\cos \varphi| \le 1/\sqrt{2}$ and $|\cos \varphi| \ge 1/\sqrt{2}$.
In the first case we have $\sin \varphi \ge 1/\sqrt{2}$ and thus
 \[ (1+|v|\sin \varphi)^{q-(N-2)} \le C \, (1+|v|)^{q-(N-2)} \]
and the end of the proof is strictly similar to above.
In the second case, we have
 \[ \frac{(r+2|v| \cos \varphi)^2}8 \ge \frac{|v|^2}{12} - \frac{r^2}{16}  \]
which implies exponential decay in $v$ for $J'$.
\end{proof}

\section{Proof of Theorem~\ref{theo:mainB} for $\var >0$}
\label{sec:Bconst}
\setcounter{equation}{0}

Let us first consider a (non locally integrable) collision kernel of the form (in $\sigma$-representation)
 \[ B_{\gamma,\alpha} (|v-v_*|,\cos \theta) = |v-v_*|^\gamma \, \sin^{-(N-1)-\alpha} (\theta/2) \]
with $\gamma \in (-N,+\infty)$ and $\alpha \in [0,2)$.

The {\em Dirichlet form} of the corresponding linearized collision operator $L_{\cal B}$ is
 \begin{multline*}
 D_{\gamma, \alpha} (g) := \langle L_{\cal B} g, g \rangle_{L^2(\R^N)} \\
          = \frac14 \, \int_{\R^N \times \R^N \times \ens{S}^{N-1}}
           \left[ \left(\frac{g(v')}{\mu(v')^{1/2}}\right)
      + \left(\frac{g(v'_*)}{\mu(v'_*)^{1/2}}\right) - \left(\frac{g(v)}{\mu(v)^{1/2}}\right)
      - \left(\frac{g(v_*)}{\mu(v_*)^{1/2}}\right) \right]^2 \\ \times
        |v-v_*|^\gamma \, \sin^{-(N-1)-\alpha} (\theta/2) \, \mu(v) \, \mu(v_*) \, dv \, dv_* \, d\sigma.
 \end{multline*}
 
 We want to prove the
 \begin{proposition}\label{prop:mainBE}
 For any $\var >0$, there is an explicit constant $C_{\gamma,\alpha,\var}>0$ such that
   \[  D_{\gamma, \alpha} (g) \ge 
       C_{\gamma,\alpha,\var} \, \int_{\R^N} [g - {\bf P}g ]^2 \, (1+|v|)^{\gamma + \alpha -\var} \, dv. \]
 \end{proposition}

\begin{proof}[Proof of Proposition~\ref{prop:mainBE}]
For any $\beta \in (0,(N-1)+\alpha)$, we introduce the following angular truncation domain
(which depends on $|v-v_*|$):
 \[ \cal{C}_{\beta} = \left\{ \sigma \in \ens{S}^{N-1} \, ; \
       \sin^{-(N-1)-\alpha} (\theta/2) \ge |v-v_*|^{\beta} \right\}. \]
One checks easily that $\cal{C}_{\beta}$ is invariant under the
pre-post collisional change of variables and the change of variable
$(v,v_*,\sigma) \to (v_*,v,-\sigma)$ (for the different classical changes of variable 
we refer to~\cite[Chapter~1, Section~4]{Vi:hb}).

Hence we have
 \begin{multline*}
 D_{\gamma, \alpha} (g) \ge \frac14 \, \int_{\R^N \times \R^N \times \cal{C}_{\beta}}
           \left[ \left(\frac{g(v')}{\mu(v')^{1/2}}\right)
      + \left(\frac{g(v'_*)}{\mu(v'_*)^{1/2}}\right) - \left(\frac{g(v)}{\mu(v)^{1/2}}\right)
      - \left(\frac{g(v_*)}{\mu(v_*)^{1/2}}\right) \right]^2 \\ \times
        |v-v_*|^{\gamma+\beta} \, \mu(v) \, \mu(v_*) \, dv \, dv_* \, d\sigma \\
   = \int_{\R^N \times \R^N \times \cal{C}_{\beta}}
           \left[ - \left(\frac{g(v')}{\mu(v')^{1/2}}\right)
      - \left(\frac{g(v'_*)}{\mu(v'_*)^{1/2}}\right) + \left(\frac{g(v)}{\mu(v)^{1/2}}\right)
      + \left(\frac{g(v_*)}{\mu(v_*)^{1/2}}\right) \right] \, \left(\frac{g(v)}{\mu(v)^{1/2}}\right) \\ \times
        |v-v_*|^{\gamma+\beta} \, \mu(v) \, \mu(v_*) \, dv \, dv_* \, d\sigma \\
   \ge \int_{\R^N \times \R^N \times \cal{C}_{\beta}} g^2 \, |v-v_*|^{\gamma+\beta} 
           \, \mu(v_*) \, dv \, dv_* \, d\sigma \\
       + \int_{\R^N \times \R^N \times \cal{C}_{\beta}} g \, g_* \, |v-v_*|^{\gamma+\beta} \,
           \mu^{1/2}(v) \mu^{1/2}(v_*) \, dv \, dv_* \, d\sigma \\
       - \int_{\R^N \times \R^N \times \cal{C}_{\beta}} g \, g' \, |v-v_*|^{\gamma+\beta} \,
           \mu^{1/2}(v_*) \mu^{1/2} (v'_*) \, dv \, dv_* \, d\sigma \\
       - \int_{\R^N \times \R^N \times \cal{C}_{\beta}} g \, g'_* \, |v-v_*|^{\gamma+\beta} \,
           \mu^{1/2}(v_*) \mu^{1/2}(v'_*) \, dv \, dv_* \, d\sigma \\ =: D_1 + D_2 + D_3 + D_4.
 \end{multline*}

In the following we shall bound $D_1$ from below, and $D_2, D_3, D_4$ from above.
\smallskip

For $D_1$ we have
 \[ D_1 \ge \int_{\R^N \times \R^N} g^2 \, |v-v_*|^{\gamma+\beta} \, \mu(v_*) \,
                  \left(\int_{\sigma \in \cal{C}_\beta} d\sigma \right) \, dv \, dv_*. \]
An easy computation yields
 \[ 
 \int_{\sigma \in \cal{C}_\beta} d\sigma = C \, \int_0 ^a \sin^{N-2} \theta \, d\theta 
 \]
with $a = 2 \arcsin \big( |v-v_*|^{-\beta/((N-1)+\alpha)} \big) \ge C \,|v-v_*|^{-\beta/((N-1)+\alpha)}$,
and thus
 \[ \int_{\sigma \in \cal{C}_\beta} d\sigma \ge C \, \int_0 ^a \theta^{N-2} \, d\theta =
    C \, a^{N-1} \ge C \, |v-v_*|^{-\frac{\beta \, (N-1)}{(N-1)+\alpha}}.  \]

Hence we deduce
 \begin{multline*}
 D_1 \ge C \, \int_{\R^N} g^2 \,
    \left( \int_{\R^N} |v-v_*|^{\gamma+\frac{\beta \alpha}{(N-1) + \alpha}} \, \mu(v_*) \, dv_* \right) \, dv \\
      \ge C \, \int_{\R^N} g^2 \, (1+|v|)^{\gamma+\frac{\beta \alpha}{(N-1) + \alpha}} \, dv .
 \end{multline*}
This completes the estimate for the first term $D_1$.  
\smallskip

Let us consider the second term $D_2$. 
  \begin{multline*}
  D_2 = \int_{\R^N \times \R^N} g \, g_* \, |v-v_*|^{\gamma+\beta} \, 
            \mu^{1/2}(v) \mu^{1/2}(v_*) \, dv \, dv_* \\ \le 
            \int_{\R^N \times \R^N} g^2 \, \mu^{1/2}(v) \mu^{1/2}(v_*) \, 
                |v-v_*|^{\gamma+\beta} \, dv \, dv_* \\ 
         \le C \, \left( \int_{\R^N} g^2 \, (1+|v|)^\gamma \, dv \right). 
  \end{multline*}
\smallskip

Let us consider the terms $D_3$ and $D_4$. We remove the truncation $\cal{C}_\beta$ by
bounding from above by the collision kernel without truncation, and we have
 \[ |D_3| + |D_4| \ \le \int_{\R^N \times \R^N} |g(v)| \, |g(v')| \, k_{\gamma + \beta}(v,v') \, dv \, dv' \]
(once the truncation is removed, the kernel is invariant under the change of variable 
$(v,v_*,\sigma) \to (v,v_*,-\sigma)$ which allows to reduce to the same term).

Hence we can write
 \begin{multline*}
 |D_3| + |D_4| \  \le C \, \left( \int_{\R^N} g^2 \, (1+|v|)^{\gamma + \beta -(N-1)} \, dv \right)^{1/2} \\ \times
            \left[ \int_{\R^N} (1+|v|)^{-(\gamma + \beta -(N-1))} \,
               \left( \int_{\R^N} k_{\gamma + \beta}(v,v') \, |g(v')| \, dv' \right)^2 \, dv \right]^{1/2} \\
     \le C \, \left( \int_{\R^N} g^2 \, (1+|v|)^{\gamma + \beta -(N-1)} \, dv \right)^{1/2} \\ \times
           \left[ \int_{\R^N} (1+|v|)^{-(\gamma + \beta -(N-1))} \,
               \left( \int_{\R^N} k_{\gamma + \beta}(v,v') \, dv' \right) \left(
                      \int_{\R^N} k_{\gamma + \beta}(v,v'')\, g(v'')^2 \, dv'' \right) \, dv \right]^{1/2}.
 \end{multline*}

Then we want to use Proposition~\ref{prop:Grad} to get
\begin{equation} \label{controlk}
(1+|v|)^{-(\gamma + \beta -(N-1))} \, 
\left( \int_{\R^N} k_{\gamma + \beta}(v,v') \, dv' \right) \le C. 
\end{equation}
It is possible as soon as $\gamma+ \beta >-1$,
{\it i.e.}, $\beta > -\gamma -1$. Since $\gamma > -N$ it is enough
that $\beta > (N-1)$. Since $\alpha >0$ by assumption, it is always possible
to pick $\beta$ such that
 \[ (N-1) < \beta < (N-1) + \alpha. \]
For this choice of $\beta$, one can apply Proposition~\ref{prop:Grad} to get~\eqref{controlk}.

Thus
 \begin{multline*}
 |D_3| + |D_4| \ \le C \, \left( \int_{\R^N} g^2 \, (1+|v|)^{\gamma + \beta -(N-1)} \, dv \right)^{1/2} \\ \times
           \left[ \int_{\R^N}
     \left( \int_{\R^N} k_{\gamma + \beta}(v,v'')\, g(v'')^2 \, dv'' \right) \, dv \right]^{1/2}  \\
     \le C \, \left( \int_{\R^N} g^2 \, (1+|v|)^{\gamma + \beta -(N-1)} \, dv \right)^{1/2} \\ \times
           \left[ \int_{\R^N}
     \left( \int_{\R^N} k_{\gamma + \beta}(v,v'') \, dv \right) \, g(v'')^2 \, dv'' \right]^{1/2}.
 \end{multline*}
Using again Proposition~\ref{prop:Grad} we have
 \[ \left( \int_{\R^N} k_{\gamma + \beta}(v,v'') \, dv \right) \le C \, (1+|v''|)^{\gamma + \beta -(N-1)} \]
which yields finally
 \[ |D_3| + |D_4| \ \le C \, \left( \int_{\R^N} g^2 \, (1+|v|)^{\gamma + \beta -(N-1)} \, dv \right). \]

But the choice $\beta \in ((N-1),(N-1) + \alpha)$ implies straightforwardly
 \[ \gamma < \gamma + \beta -(N-1) < \gamma + \frac{\beta \alpha}{(N-1) + \alpha}. \]
Hence by trivial interpolation we get for any $\delta >0$ there is an explicit 
constant $C_\delta >0$ such that 
 \[ |D_3| + |D_4| \ \le C_\delta \, \left( \int_{\R^N} g^2 \, (1+|v|)^{\gamma} \, dv \right)
           + \delta \, \left( \int_{\R^N} g^2 \, (1+|v|)^{\gamma +\frac{\beta \alpha}{(N-1) + \alpha}} \, dv \right). \]
Then by taking $\delta$ smaller than the constant in the bound from below for $D_1$ we deduce that
 \[ D_{\gamma, \alpha} (g) \ge 
      C_+ \, \left( \int_{\R^N} g^2 \, (1+|v|)^{\gamma +\frac{\beta \alpha}{(N-1) + \alpha}} \, dv \right)
             - K_- \, \left( \int_{\R^N} g^2 \, (1+|v|)^\gamma \, dv \right). \]

To conclude the proof we finally use the coercivity estimates in~\cite{Mcoerc}. In this paper
it is proved that (under our assumption on the collision kernel)
 \[ D_{\gamma, \alpha} (g) \ge C_0 \, \left( \int_{\R^N} \big[g - {\bf P}g \big]^2 \, (1+|v|)^\gamma \, dv \right) \]
for some explicit constant $C_0 >0$.

By combining these two last inequalities, one deduces that
 \[ D_{\gamma, \alpha} (g) \ge 
    C_1 \, \left( \int_{\R^N} \big[g - {\bf P}g \big]^2 \, (1+|v|)^{\gamma+\frac{\beta \alpha}{(N-1) + \alpha}} \, dv \right) \]
for some explicit constant $C_1 >0$ (depending on $\beta$).
Since $\beta$ can be taken as close as wanted to $(N-1) + \alpha$, the weight
exponent $\gamma+ \beta \alpha/((N-1) + \alpha)$ can be taken as close as wanted to $\gamma + \alpha$.
This concludes the proof.
\end{proof}

\begin{proof}[Proof of the first point in Theorem~\ref{theo:mainB}]
To conclude the proof of the first point in Theorem~\ref{theo:mainB}, it is enough to remark that the proof 
of Proposition~\ref{prop:mainBE} can be modified easily in order to start from a  
collision kernel $B$ such that 
$$
B \ge K \, B_{\gamma,\alpha} \, {\bf 1}_{\theta \in [0,\theta_0]}
$$
for some constants $K>0$, $\theta_0 \in (0,\pi]$ and 
 \[ B_{\gamma,\alpha} (|v-v_*|,\cos \theta) = 
     |v-v_*|^\gamma \, \sin^{-(N-1)-\alpha} (\theta/2), \quad \gamma \in (-N,+\infty), \ \alpha \in [0,2), \]
defined as in the beginning of this section. The assumptions of Theorem~\ref{theo:mainB} 
imply such a control.

Indeed one first reduces to the collision kernel 
$K \, B_{\gamma,\alpha} \, {\bf 1}_{\theta \in [0,\theta_0]}$ by monotonicity of the 
Dirichlet form. 
Then the bound from above on $D_2$, $D_3$, $D_4$ are unchanged, and 
the bound from below on $D_1$ is still valid since one gets straightforwardly 
$$
\int_{\sigma \in \cal{C}_\beta}  {\bf 1}_{\theta \in [0,\theta_0]} \, d\sigma 
\ge \min \Big\{ C_1 \, |v-v_*|^{- \frac{\beta (N-1)}{(N-1) + \alpha}}; C_2 \Big\} 
$$
for some constants $C_1,C_2>0$. Finally the coercivity estimates of~\cite{Mcoerc} 
used in the proof are also still valid for a collision kernel 
$K \, B_{\gamma,\alpha} \, {\bf 1}_{\theta \in [0,\theta_0]}$.
\end{proof}

\section{Proof of Theorem~\ref{theo:mainB} for $\var=0$} 
\label{sec:Bnc}
\setcounter{equation}{0}

We start again from some collision kernel which satisfies 
$$
B \ge K \, B_{\gamma,\alpha} \, {\bf 1}_{\theta \in [0,\theta_0]}
$$
for some constants $K>0$, $\theta_0 \in (0,\pi]$ and 
 \[ B_{\gamma,\alpha} (|v-v_*|,\cos \theta) = 
    |v-v_*|^\gamma \, \sin^{-(N-1)-\alpha} (\theta/2), \quad \gamma \in (-N,+\infty), \ \alpha \in [0,2). \]

We shall prove the
 \begin{proposition}\label{prop:border}
 There is some constant $C_{B,0} >0$ (obtained in our proof by compactness argument) such that
   \[  D_{B} (g) \ge C_{B,0} \, \int_{\R^N} \big[g - {\bf P}g \big]^2 \, (1+|v|)^{\gamma + \alpha} \, dv. \]
 \end{proposition}
 
The second point in Theorem~\ref{theo:mainB} follows 
immediately from Proposition~\ref{prop:border}.

\begin{proof}[Proof of Proposition~\ref{prop:border}]
We introduce the following angular truncation domain (which depends on $|v-v_*|$):
 \[ \cal{\bar C} = \left\{ \sigma \in \ens{S}^{N-1} \, ; \ |v-v'| \le 1 \right\}. \]
One checks easily that $\cal{\bar C}$ is invariant under the
pre-post collisional change of variables and the change of variable
$(v,v_*,\sigma) \to (v_*,v,-\sigma)$. Remark that this truncation domain 
$\cal{\bar C}$ corresponds to the limit case $\beta = (N-1) + \alpha$ in the 
truncation domain $\cal{C}_\beta$ previously introduced. 

Hence we have
 \begin{multline*}
 D_{\gamma, \alpha} (g) \ge \frac{K}4 \, \int_{\R^N \times \R^N \times \cal{\bar C}}
           \left[ \left(\frac{g(v')}{\mu(v')^{1/2}}\right)
      + \left(\frac{g(v'_*)}{\mu(v'_*)^{1/2}}\right) - \left(\frac{g(v)}{\mu(v)^{1/2}}\right)
      - \left(\frac{g(v_*)}{\mu(v_*)^{1/2}}\right) \right]^2 \\ \times
        |v-v_*|^{\gamma+\alpha+(N-1)} \, {\bf 1}_{\theta \in [0,\theta_0]} \, \mu(v) \, \mu(v_*) \, dv \, dv_* \, d\sigma \\
   = K \, \int_{\R^N \times \R^N \times \cal{\bar C}}
           \left[ - \left(\frac{g(v')}{\mu(v')^{1/2}}\right)
      - \left(\frac{g(v'_*)}{\mu(v'_*)^{1/2}}\right) + \left(\frac{g(v)}{\mu(v)^{1/2}}\right)
      + \left(\frac{g(v_*)}{\mu(v_*)^{1/2}}\right) \right] \, \left(\frac{g(v)}{\mu(v)^{1/2}}\right) \\ \times
        |v-v_*|^{\gamma+\alpha+(N-1)} \, {\bf 1}_{\theta \in [0,\theta_0]} \, \mu(v) \, \mu(v_*) \, dv \, dv_* \, d\sigma \\
   = K \, \langle \hat L g, g \rangle
 \end{multline*}
with a fictious self-adjoint operator $\hat L$ on $L^2$ defined by 
\begin{multline*} 
\hat L g (v) = \int_{\R^N \times \ens{S}^{N-1}}
           \Bigg[ - \left(\frac{g(v')}{\mu(v')^{1/2}}\right)
      - \left(\frac{g(v'_*)}{\mu(v'_*)^{1/2}}\right) + \left(\frac{g(v)}{\mu(v)^{1/2}}\right)
      + \left(\frac{g(v_*)}{\mu(v_*)^{1/2}}\right) \Bigg] \\ 
        |v-v_*|^{\gamma+\alpha+(N-1)} \, {\bf 1}_{|v-v'| \le 1} \, {\bf 1}_{\theta \in [0,\theta_0]} \, 
        \mu^{1/2} (v) \, \mu(v_*) \, dv_* \, d\sigma.
\end{multline*}

Using the same decomposition as in Section~\ref{sec:tech}, this operator can be split as 
$$
\hat L = \hat \nu - \hat K^+ + \hat K^c 
$$
where the multiplicative function $\hat \nu$ is 
$$
\hat \nu (v) = \int_{\R^N \times \ens{S}^{N-1}} 
|v-v_*|^{\gamma+\alpha+(N-1)} \, {\bf 1}_{|v-v'| \le 1} \, {\bf 1}_{\theta \in [0,\theta_0]} \, \mu(v_*) \, dv_* \, d\sigma
$$
which satisfies by similar computations as above 
$$
\hat \nu(v) \ge C \, (1 + |v|)^{\gamma + \alpha}, \quad C> 0.
$$

Then we shall show that the remaining terms satisfies some compactness property. 
\smallskip

Let us assume first for the sake of clarity that $\gamma + \alpha =0$. Then the multiplication 
function $\hat \nu$ is bounded from below by a positive constant $\hat \nu_0 >0$, 
and it is straightforward that $\hat K^c$ is a Hilbert-Schmidt operator (we leave these computations 
to the reader). 
Let us show that the part $\hat K^+$ can be written as a limit of 
Hilbert-Schmidt operators (hence showing that it is compact in $L^2$). 
The kernel of $\hat K^+$ is by inspection 
\begin{multline*}
\hat k := k_{(N-1)} (v,v') \, {\bf 1}_{|v-v'| \le 1} \, {\bf 1}_{\theta \in [0,\theta_0]} \\
=  \frac{2^N}{|v'-v| \, (2 \pi)^{3/2}} \,
             \exp\left\{ -\frac{|v'-v|^2}8 -\frac{|v'-v + 2 (v \cdot \omega) \omega|^2}8 \right\} \\ \times
                 \left( \int_{\omega^\bot} |v'-v + z| \, 
                 \exp\left\{ - \frac{|z + \left(v - (v \cdot \omega) \omega \right)|^2}2 \right\} \, dz \right) 
                 \, {\bf 1}_{|v-v'| \le 1} \, {\bf 1}_{\theta \in [0,\theta_0]}
\end{multline*}
(where $k_q$, $q \in (-N,1]$ was the kernel computed in Section~\ref{sec:tech}). 
We approximate this kernel (and correspondingly the operator $\hat K^+$) as follows: 
$$ 
\hat k = \hat k_\var ^c + \hat k_\var ^r
$$
with 
$$
\hat k_\var ^c = \left({\bf 1}_{|v-v'| \ge \var} 
\times  {\bf 1}_{\left|\frac{v}{|v|} \cdot \frac{(v-v')}{|v-v'|}\right| \ge \var} \right) \, \hat k 
$$
and 
$$
\hat k_\var ^r = \hat k  - \hat k_\var ^c.
$$
By similar straightforward computations as in the proof of Proposition~\ref{prop:Grad} we get that 
$\hat k_\var ^r$ is symmetric in $v,v'$ and 
$$
\sup_{v \in \R^N} \int_{\R^N} |\hat k_\var ^r| \, dv'  \xrightarrow[]{\var \to 0} 0 
$$
and therefore $\hat K^{+,c} _\var \to \hat K^+$ in $L^2$ as $\var \to 0$. Hence it is enough to show that 
$\hat K^{+,c}$ is compact. But the kernel $\hat k_\var ^c$ satisfies (using the same representation 
as in Proposition~\ref{prop:Grad})
\begin{multline*} 
\int_{\R^N \times \R^N} \big( \hat k_\var ^c \big)^2 \, dv \, dv' \le 
C \, \int_{v \in \R^N} 
\Bigg( \int_\var ^1 r^{N-3} \, (1+r)^2 \, (1+|v| \, \sin \theta)^2 \,
    e^{-\frac{r^2}4} \\ \int_0 ^\pi e^{-\frac{(r+2|v| \cos \theta)^2}4} \, \sin^{N-2} \theta \,  {\bf 1}_{|\cos \theta | \ge \var} \, d\theta 
\Bigg) \\ 
\le C \, \int_{v \in \R^N}  (1+|v|)^2 \, e^{- \var^2 \, |v|^2} \, dv < +\infty. 
\end{multline*} 
Therefore $\hat K^{+,c}$ is a Hilbert-Schmidt operator and the result is proved. 
By applying Weyl's theorem (exactly as in~\cite{Grad63,CIP}), we deduce that 
$\hat L$ and $\hat \nu$ have the same essential spectrum, which is included 
in $[\hat \nu_0, +\infty)$, and therefore, since $\hat L \ge 0$, that $0$ is an isolated 
eigenvalue, which concludes the proof.
\smallskip

When $\gamma + \alpha$ is different from $0$, one considers (in the spirit of~\cite{GoPo86}) 
the following symmetric weighted modification of $\hat L$: 
$$
\tilde L = (1 + |\cdot|)^{-(\gamma + \alpha)/2} \, 
             \hat L \, \Big( (1 + |\cdot|)^{-(\gamma + \alpha)/2} \cdot \Big)
$$
and the corresponding splitting $\tilde L = \tilde \nu -\tilde K^+ + \tilde K^c$. 
Then $\tilde \nu$ is strictly positive uniformly (and bounded from above) and it is 
straightforward again that $\tilde K^c$ is a Hilbert-Schmidt operator. Let us focus 
on the term $\tilde K^+$. Its kernel is 
$$
(1 + |v|)^{-(\gamma + \alpha)/2} \, k_{\gamma + \alpha + (N-1)} (v,v') 
\, (1 + |v'|)^{-(\gamma + \alpha)/2} \, {\bf 1}_{|v-v'| \le 1} \, {\bf 1}_{\theta \in [0,\theta_0]}
$$
and similar computations as above show again that $\tilde K^+$ can be written as a limit 
of Hilbert-Schmidt operators (remark 
that thanks to the truncation $|v-v'| \le 1$, weights on $v$ and $v'$ can be interchanged up to a constant). 

We then conclude by applying Weyl's Theorem to $\tilde L$ (as in~\cite{GoPo86}). We deduce 
thus that $0$ is an isolated eigenvalue of $\tilde L \ge 0$, and therefore we obtain the existence 
of the inequality in Proposition~\ref{prop:border}.  
\end{proof}




\section{Proof of Theorem~\ref{theo:mainL}}
\label{sec:L}
\setcounter{equation}{0}

In \cite{Guo:Land} it was shown {\it via} compactness arguments that there is 
a constant $C^1 _\gamma >0$ such that 
\begin{equation}
\langle L^{\mathcal{L}}g,g\rangle \ge C^1 _\gamma \, \big\| [g - {\bf P}g ] \big\|_\sigma ^2.  \label{vlower}
\end{equation}
where $\|\cdot \|_\sigma$ is the following anisotropic norm: 
$$
\|g\|_{\sigma}^2 := \int_{\mathbb{R}^N}\Big(  (1+|v|)^{\gamma} \, \big|\Pi_v\nabla_v g\big|^2
+ (1+|v|)^{\gamma +2} \, \big|[\mbox{Id} -\Pi_v]\nabla_v g\big|^2+(1+|v|)^{\gamma+2}\, g^2 \Big) \, dv
$$
with 
$$
\Pi_v \nabla_v g= \left(\frac{v}{|v|}\cdot \nabla_v g\right) \frac{v}{|v|}.
$$

In particular this lower bound implies a spectral gap for the linearized Landau collision operator 
as soon as $\gamma\ge -2$.  
\smallskip

On the other hand, in~\cite{Mcoerc} the first author derived the following explicit coercivity estimate:
\begin{equation}
\langle L^{\mathcal{L}}[\mu^{1/2}h],[\mu^{1/2}h]\rangle \ge C_\gamma^2 \, 
\int_{\mathbb{R}^N}(1+|v|^2)^{\gamma/2}\Big(  |\nabla_v h|^2+h^2 \Big)\mu(v) \, dv
= C_\gamma ^2 \, \| h \|_{H^1_\gamma(\mu)}.
\label{mest}
\end{equation}
This holds for all $h= h-  {\bf {\bf \hat P}} h$, where ${\bf \hat P}$ is the orthogonal projection in $L^2(\mu)$ 
given by
$$
{\bf \hat P} h= a+b\cdot v+c|v|^2
$$
with $a,c\in\mathbb{R}$ and $b\in \mathbb{R}^N$. 
This estimate originates from the alternate (but equivalent) linearization
$$
f=\mu(1 +h).
$$


\begin{proof}[Proof of Theorem~\ref{theo:mainL}]
Letting $h=\mu^{-1/2}g$, we first translate the estimate of~\cite{Mcoerc} into a coercivity estimate for $g$:  
the r.h.s. of \eqref{mest} is given by
\begin{equation*}
\begin{split}
\| \mu^{-1/2}g \|_{H^1_\gamma(\mu)}
&=
\int_{\mathbb{R}^N}(1+|v|^2)^{\gamma/2} \, \left(  \left|\nabla_v g+\frac{1}{2}v g\right|^2+g^2 \right) \, dv
\\
&=
\int_{\mathbb{R}^N}(1+|v|^2)^{\gamma/2} \, \left(  \left|\nabla_v g\right|^2+\frac{1}{4} \, |v|^2 \, |g|^2 
+ (v\cdot \nabla_v g) \, g+g^2 \right) \, dv.
\end{split}
\end{equation*}
We focus our attention on the term without a definite sign:
\begin{equation*}
\begin{split}
\int_{\mathbb{R}^N}(1+|v|^2)^{\gamma/2}(v\cdot\nabla_v g)g \, dv
&=
\frac{1}{2}\int_{\mathbb{R}^N}(1+|v|^2)^{\gamma/2}v\cdot\nabla_v (g^2) \, dv
\\
&=
-\frac{1}{2}\int_{\mathbb{R}^N}\nabla\cdot \left( v(1+|v|^2)^{\gamma/2}\right) g^2 \, dv.
\end{split}
\end{equation*}
Now we look at the derivative of the polynomial 
$$
\nabla\cdot \left( v(1+|v|^2)^{\gamma/2}\right)
=
(1+|v|^2)^{\gamma/2}\left(N+\frac{\gamma \, |v|^2}{1+|v|^2}\right).
$$
Further notice that 
$$
\left|N+\frac{\gamma \, |v|^2}{1+|v|^2}\right|\le N+|\gamma|.  
$$
Plug in these last few computations to obtain
\begin{multline*}
\| \mu^{-1/2}g \|_{H^1_\gamma(\mu)}
\ge
\frac{1}{4}\int_{\mathbb{R}^N} \left( (1+|v|^2)^{\frac{\gamma}{2}} \, \left|\nabla_v g\right|^2 
+ (1+|v|^2)^{\frac{\gamma+2}{2}} \, g^2 \right) \, dv \\
-\frac{1}{2}(N+|\gamma|)\int_{\mathbb{R}^N}(1+|v|^2)^{\frac{\gamma}{2}}g^2 \, dv.
\end{multline*}
On the other hand 
$$
\| \mu^{-1/2}g \|_{H^1_\gamma(\mu)}
\ge
\int_{\mathbb{R}^N}(1+|v|^2)^{\frac{\gamma}{2}} \, g^2 \, dv.
$$
Combining the two last inequalities we thus conclude 
$$
\| \mu^{-1/2}g \|_{H^1_\gamma(\mu)}
\ge
C \, \int_{\mathbb{R}^N}\Big[ (1+|v|^2)^{\frac{\gamma}{2}} \, \left|\nabla_v g\right|^2+(1+|v|^2)^{\frac{\gamma+2}{2}} \, g^2 \Big] \, dv.
$$

Combine this with \eqref{mest} to obtain an explicitly computable constant $C>0$ 
such that for any $g=g - {\bf P}g$, we have 
\begin{equation}
\langle L^{\mathcal{L}}g,g\rangle \ge C
\int_{\mathbb{R}^N} \left( (1+|v|^2)^{\frac{\gamma}{2}}  \left|\nabla_v g\right|^2+(1+|v|^2)^{\frac{\gamma+2}{2}} \, g^2 \right) \, dv.
\label{mestg}
\end{equation}
While the power is the same as in \eqref{vlower} for the term with no derivative, 
the derivative term still has a better power in \eqref{vlower}.

Guo \cite{Guo:Land} computes that 
$$
\langle L^{\mathcal{L}}g,g\rangle=\|g\|_\sigma ^2-\langle \partial _i\sigma ^ig ,g\rangle
-\langle Kg,g\rangle. 
$$
Here we will not give precise definitions of $\partial _i\sigma ^i$ and $K$, we will only use an estimate with explicitly computable constants for these terms to get an explicit lower bound in the $\sigma$ norm.  

The \cite[Lemma 5]{Guo:Land} implies that for any $m>1,$ there is an (explicit) 
$0<C(m)<\infty$ such that 
\begin{equation*}
|\langle \partial _i\sigma ^ig,g\rangle |+|\langle
Kg_1,g_2\rangle |  
\le 
\frac{1}{m}\left\| g\right\| _{\sigma} ^2
+
C(m) \, \left\| {\bf 1}_{\{|\cdot| \le C(m)\}} g \right\|^2 _{L^2}. 
\end{equation*}
We deduce 
\begin{equation*}
|\langle \partial _i\sigma ^ig,g\rangle |+|\langle Kg_1,g_2\rangle |  
\le 
\frac{1}{m}\left\| g\right\| _{\sigma} ^2
+ C'(m) \, \langle L^{\mathcal{L}}g,g\rangle  
\end{equation*}
for another explicit constant $C'(m)>0$ thanks to~\eqref{mestg}. Taking 
for instance $m=2$ we deduce
$$
\left\| g\right\| _{\sigma} ^2 = \langle L^{\mathcal{L}}g,g\rangle + \langle \partial _i\sigma ^ig ,g\rangle
+ \langle Kg,g\rangle \le (1+C'(2)) \, \langle L^{\mathcal{L}}g,g\rangle + \frac12 \, \left\| g\right\| _{\sigma} ^2,
$$
which concludes the proof.
\end{proof}

%

\bigskip
\noindent
{\bf{Acknowledgments:}} This work was initiated while the first author was visiting Brown University, and 
he wishes to thank its Department of Applied Mathematics, in particular Yan Guo and Robert Strain for 
their hospitality. 
The second author, while finishing this work, was partially supported by a National Science 
Foundation Mathematical Sciences Postdoctoral Research Fellowship.

\begin{flushleft} \signcm \end{flushleft}
\begin{flushleft} \signrs \end{flushleft}

\smallskip


\begin{thebibliography}{100}

\bibitem{AlexTTSP}
{\sc Alexandre, R.} 
\newblock Remarks on 3D Boltzmann linear equation without cutoff. 
\newblock {\it Transport Theory Statist. Phys. 28}, 5 (1999), 433--473.

\bibitem{AlexCRASrenorm}
{\sc Alexandre, R.} 
\newblock Une d\'efinition des solutions renormalis\'ees pour l'\'equation de Boltzmann sans troncature angulaire. 
\newblock {\it C. R. Acad. Sci. Paris S\'er. I Math. 328}, 11 (1999), 987--991.

\bibitem{AlexM2AN}
{\sc Alexandre, R.} 
\newblock Around 3D Boltzmann non linear operator without angular cutoff, a new formulation.
\newblock {Math. Model. Numer. Anal. 34}, 3 (2000), 575--590.

\bibitem{AlexM3AS}
{\sc Alexandre, R. and El Safadi, M.}
\newblock Littlewood-Paley theory and regularity issues in Boltzmann homogeneous equations. I. Non-cutoff case and Maxwellian molecules. 
\newblock {\it Math. Models Methods Appl. Sci. 15}, 6 (2005), 907--920.

\bibitem{ADVW}
{\sc Alexandre, R., Desvillettes, L., Villani, C., Wennberg, B.}
\newblock Entropy dissipation and long-range interactions. 
\newblock {\it Arch. Rational Mech. Anal. 152} (2000), 327--355.

\bibitem{AlViCPAM} 
{\sc Alexandre, R. and Villani, C.}
\newblock On the Boltzmann equation for long-range interactions.
\newblock {\it Comm. Pure Appl. Math. 55}, 1 (2002), 30--70.

\bibitem{AlVi:04}
{\sc Alexandre, R., Villani, C.},
\newblock On the {L}andau approximation in plasma physics. 
\newblock {\it Ann. Inst. H. Poincar\'e Anal. Non Lin\'eaire 21}, 1 (2004), 61--95. 


\bibitem{ArkInf1} 
{\sc Arkeryd, L.} 
\newblock Intermolecular forces of infinite range and the Boltzmann equation.
\newblock {\it Arch. Rational Mech. Anal. 77}, 1 (1981), 11--21.

\bibitem{ArkInf2} 
{\sc Arkeryd, L.} 
\newblock Asymptotic behaviour of the Boltzmann equation with infinite range forces.
\newblock {\it Comm. Math. Phys. 86}, 4 (1982), 475--484.


\bibitem{ArBu:91}
{\sc Arsen$'$ev, A. A., Buryak, O. E.}
\newblock On a connection between the solution of the {B}oltzmann
equation and the solution of the {L}andau-{F}okker-{P}lanck equation. 
\newblock {\it Math. USSR Sbornik 69}, 2 (1991), 465--478. 

\bibitem{BaMo}
{\sc Baranger, C. and Mouhot, C.}
\newblock  Explicit spectral gap estimates for the
linearized {B}oltzmann and {L}andau operators with hard potentials.
\newblock {\it Rev. Matem. Iberoam. 21} (2005), 819--841. 

\bibitem{Boby88}
{\sc Bobyl{\"e}v, A. V.},
The theory of the nonlinear spatially uniform {B}oltzmann
equation for {M}axwell molecules.
\newblock {\it Mathematical physics reviews, Vol.\ 7} (1988),
Soviet Sci. Rev. Sect. C Math. Phys. Rev., 111--233.


\bibitem{BobyCerc:conjCer:99}
{\sc Bobylev, A. V. and Cercignani, C.}
\newblock On the rate of entropy production for the {B}oltzmann equation.
\newblock {\it J. Statist. Phys. 94} (1999), 603--618.




\bibitem{Cafl80}
{\sc Caflisch, R. E.}
\newblock The {B}oltzmann equation with a soft potential. {I}. {L}inear,
spatially-homogeneous.
\newblock {\it Comm. Math. Phys. 74} (1980), 71--95.


\bibitem{Carl:57}
{\sc Carleman, T.}
\newblock Probl\`emes math\'ematiques dans la th\'eorie cin\'etique des gaz. 
\newblock Almqvist and Wiksells Boktryckeri Ab, Uppsala 1957.





\bibitem{Ce67}
{\sc Cercignani, C.}
\newblock On Boltzmann Equation with Cutoff Potentials.
\newblock {\it Phys. Fluids 10} (1967), 2097--2104. 

\bibitem{Ce88}
{\sc Cercignani, C.}
\newblock The Boltzmann equation and its applications.
\newblock Applied Mathematical Sciences, 67. Springer-Verlag, New York, 1988.

\bibitem{CzPa} 
{\sc Czechowski, Z. and Palczewski, A.} 
\newblock Spectrum of the Boltzmann collision operator for radial cut-off potentials. 
\newblock {\it Bull. Acad. Polon. Sci. S\'er. Sci. Tech. 28}, 9--10 (1980), 387--396.

\bibitem{CIP}
{\sc Cercignani, C. and Illner, R. and Pulvirenti, M.}
\newblock The mathematical theory of dilute gases.
\newblock Springer-Verlag, New York, 1994.

\bibitem{DeLe97}
{\sc Degond, P. and Lemou, M.}
\newblock  Dispersion relations for the linearized {F}okker-{P}lanck equation.
\newblock {\it Arch. Rational Mech. Anal.} {\bf 138} (1997), no.~2, 137--167.

\bibitem{DeLu:92}
{\sc Degond, P. and Lucquin-Desreux, B.},
\newblock The {F}okker-{P}lanck asymptotics of the {B}oltzmann collision
operator in the {C}oulomb case. 
\newblock {\it Math. Models Methods Appl. Sci. 2}, 2 (1992), 167--182. 

\bibitem{Desv:asBE:92}
{\sc Desvillettes, L.},
\newblock On asymptotics of the {B}oltzmann equation when the collisions
become grazing. 
\newblock {\it Transport Theory Statist. Phys. 21}, 3 (1992), 259--276. 

\bibitem{Desv:reg:95} 
{\sc Desvillettes, L.} 
\newblock About the regularizing properties of the non-cut-off Kac equation.
\newblock {\it Comm. Math. Phys. 168}, 2 (1995), 417--440.

\bibitem{Desv:reg:96}
{\sc Desvillettes, L.} 
\newblock Regularization for the non-cutoff $2$D radially symmetric Boltzmann equation with a velocity dependent cross section.  
\newblock Proceedings of the Second International Workshop on Nonlinear Kinetic Theories and Mathematical Aspects of Hyperbolic Systems (Sanremo, 1994). 
{\it Transport Theory Statist. Phys. 25}, 3--5 (1996), 383--394.

\bibitem{Desv:reg:97}
{\sc Desvillettes, L.} 
\newblock Regularization properties of the $2$-dimensional non-radially symmetric non-cutoff spatially homogeneous Boltzmann equation for Maxwellian molecules. 
\newblock {\it Transport Theory Statist. Phys. 26}, 3 (1997), 341--357.

\bibitem{DeVi:L2:00}
{\sc Desvillettes, L. and Villani, C.} 
\newblock On the spatially homogeneous Landau equation for hard potentials. II. $H$-theorem and applications. 
\newblock {\it Comm. Partial Differential Equations 25}, 1--2 (2000), 261--298.

\bibitem{DiLiBE} 
{\sc DiPerna, R. J. and Lions, P.-L.}
\newblock On the Cauchy problem for Boltzmann equations: global existence and weak stability.
\newblock {\it Ann. of Math. (2) 130}, 2 (1989), 321--366.


\bibitem{ElmInf}
{\sc Elmroth, T.}
\newblock Global boundedness of moments of solutions of the Boltzmann equation for forces of infinite range.
\newblock {\it Arch. Rational Mech. Anal. 82}, 1 (1983), 1--12.

\bibitem{GoPo86}
{\sc Golse, F. and Poupaud, F.}
\newblock Un r\'esultat de compacit\'e pour l'\'equation de {B}oltzmann
avec potentiel mou. {A}pplication au probl\`eme de demi-espace.
\newblock {\it C. R. Acad. Sci. Paris S\'er. I Math. 303} (1986), 585--586.

\bibitem{GoNewClass}
{\sc Goudon, T.} 
\newblock On Boltzmann equations and Fokker-Planck asymptotics: influence of grazing collisions. 
\newblock {\it J. Statist. Phys. 89}, 3--4 (1997), 751--776.

\bibitem{Grad58}
{\sc Grad, H.}
\newblock Principles of the kinetic theory of gases.
\newblock In {\it Fl\"ugge's Handbuch des Physik}, vol. XII, Springer-Verlag
(1958), pp. 205--294.

\bibitem{Grad63}
{\sc Grad, H.}
\newblock Asymptotic theory of the {B}oltzmann equation. {II}.
\newblock {\it Rarefied Gas Dynamics (Proc. 3rd Internat. Sympos., Palais de
l'UNESCO, Paris, 1962)}, Vol. I (1963), pp. 26--59.

\bibitem{Guo:VPB}
{\sc Guo, Y.}
\newblock The Vlasov-Poisson-Boltzmann system near Maxwellians. 
\newblock {\it Comm. Pure Appl. Math. 55} (2002), 1104--1135. 

\bibitem{Guo:Land}
{\sc Guo, Y.}
\newblock The Landau equation in a periodic box.
\newblock {\it Comm. Math. Phys. 231} (2002), 391--434.

\bibitem{Guo:VMB}
{\sc Guo, Y.} 
\newblock The Vlasov-Maxwell-Boltzmann system near Maxwellians. 
\newblock {\it Invent. Math. 153} (2003), 593--630. 

\bibitem{Guo:BEsoft}
{\sc Guo, Y.} 
\newblock Classical solutions to the Boltzmann equation for molecules with an angular cutoff. 
\newblock {\it Arch. Rational Mech. Anal. 169} (2003), 305--353. 

\bibitem{Guo:whole} 
{\sc Guo, Y.} 
\newblock The Boltzmann equation in the whole space. 
\newblock {\it Indiana Univ. Math. J. 53} (2004), 1081--1094. 


\bibitem{Hilb:EB:12}
{\sc Hilbert, D.}
\newblock Grundz{\"u}ge einer {A}llgemeinen {T}heorie der {L}inearen {I}ntegralgleichungen. 
\newblock {\it Math. Ann. 72}, (1912), Chelsea {P}ubl., {N}ew {Y}ork, (1953). 

\bibitem{Hinton}
{\sc Hinton, F. L.}
\newblock Collisional transport in plasma. 
\newblock {\it Handbook of Plasma Physics, Vol. 1}, North-Holland, Amsterdam, 1983.


\bibitem{Klau77}
{\sc Klaus, M.}
\newblock {B}oltzmann collision operator without cut-off.
\newblock {\it Helv. Phys. Acta 50} (1977), 893--903.



\bibitem{Lions98} 
{\sc Lions, P.-L.} 
\newblock R\'egularit\'e et compacit\'e pour des noyaux de collision de Boltzmann sans troncature angulaire.
\newblock {\it C. R. Acad. Sci. Paris S\'er. I Math. 326}, 1 (1998), 37--41.

\bibitem{Land36} 
{\sc Landau, L. D.} 
\newblock Die kinetische Gleichung f\"ur den Fall Coulombscher Wechselwirkung. 
\newblock {\it Phys. Z. Sowjet. 10} (1936), 154. 
\newblock Translation: The transport equation in the case of Coulomb interactions. 
\newblock In D. ter Haar, ed., {\em Collected papers of L. D. Landau}, pp. 163--170. 
Pergamon Press, Oxford, 1981. 




\bibitem{Mrate} 
{\sc Mouhot, C.} 
\newblock Rate of convergence to equilibrium for the spatially homogeneous Boltzmann equation with hard potentials. 
\newblock {\it Comm. Math. Phys. 261} (2006), 629--672.

\bibitem{Mcoerc}
{\sc Mouhot, C.}
\newblock Coercivity estimates for the Boltzmann and Landau operators.
\newblock To appear in {\it Comm. Partial Diff. Equations}.

\bibitem{MoNe} 
{\sc Mouhot, C. and Neumann, L.} 
\newblock Quantitative perturbative study of convergence to equilibrium for 
collisional kinetic models in the torus. 
\newblock {\it Nonlinearity 19} (2006), 969--998.


\bibitem{Pao74}
{\sc Pao, Y. P.}
\newblock {B}oltzmann collision operator with inverse-power intermolecular potentials. {I}, {II}.
\newblock {\it Comm. Pure Appl. Math. 27} (1974), 407--428, 559--581.

\bibitem{StGu:rel}
{\sc Strain, R. M., Guo, Y.} 
\newblock Stability of the relativistic Maxwellian in a collisional plasma. 
\newblock {\it Comm. Math. Phys. 251} (2004), 263--320.

\bibitem{StGu:almostexp}
{\sc Strain, R. M., Guo, Y.}
\newblock Almost Exponential Decay Near Maxwellian.
\newblock {\it Comm. Partial Diff. Equations 31} (2006), 417--429. 

\bibitem{StGu:expsoft} 
{\sc Strain, R. M., Guo, Y.} 
\newblock Exponential Decay for Soft Potentials Near Maxwellian. 
\newblock Accepted for publication in {\it Arch. Rational Mech. Anal.}. 

\bibitem{Strain} 
{\sc Strain, R. M.} 
\newblock The Vlasov-Maxwell-Boltzmann System in the Whole Space. 
\newblock Accepted for publication in {\it Comm. Math. Phys.}



\bibitem{UkAs} 
{\sc Ukai, S. and Asano, K.} 
\newblock On the Cauchy problem of the Boltzmann equation with a soft potential.
\newblock {\it Publ. Res. Inst. Math. Sci. 18}, 2 (1982), 477--519 (57--99).

\bibitem{Vill:nocu:98}
{\sc Villani, C.}
\newblock On a new class of weak solutions to the spatially homogeneous
{B}oltzmann and {L}andau equations. 
\newblock {\it Arch. Rational Mech. Anal. 143}, 3 (1998), 273--307. 

\bibitem{Vi:reg:98}
{\sc Villani, C.}
\newblock Regularity estimates via the entropy dissipation for the spatially 
homogeneous Boltzmann equation without cut-off.  
\newblock {\it Rev. Mat. Iberoamericana 15}, 2 (1999), 335--352.

\bibitem{Vi:habil}
{\sc Villani, C.}
\newblock Contribution \`a l'\'etude math\'ematique des collisions en th\'eorie cin\'etique.
\newblock Univ. Paris-Dauphine France (2000).

\bibitem{Vi:hb}
{\sc Villani, C.}
\newblock A review of mathematical topics in collisional kinetic theory.
\newblock {\it Handbook of mathematical fluid dynamics, Vol. I}, 71--305, North-Holland, Amsterdam, 2002.

\bibitem{Vi03}
{\sc Villani, C.}
\newblock Cercignani's conjecture is sometimes true and always almost true.
\newblock {\it Comm. Math. Phys. 243} (2003), 455--490.

\bibitem{WCUh52}
{\sc Wang Chang, C. S. and Uhlenbeck, G. E.}
\newblock On the propagation of Sound in Monoatomic Gases. 
\newblock Univ. of Michigan Press, Ann Arbor, Michigan 
(quoted in~\cite{Grad63}).  

\bibitem{WCUh70}
{\sc Wang Chang, C. S. and Uhlenbeck, G. E. and de Boer, J.}
\newblock Studies in Statistical Mechanics, Vol. V.
\newblock North-Holland, Amsterdam, 1970.



\end{thebibliography}
\end{document}